\documentclass{article}
\usepackage[utf8]{inputenc}
\usepackage{amsmath,amssymb}
\usepackage{graphicx}
\usepackage{authblk}
\usepackage{hyperref}
\usepackage[version=4]{mhchem}
\newcommand{\no}{\ce{NO} }
\newcommand{\noo}{\ce{NO2} }
\renewcommand{\vec}{\boldsymbol}

\title{Modelling, Simulation and Parameter Identification of Active Pollution Reduction with Photovoltaic Asphalt}

\author[2]{Jens Kruschwitz}
\author[1]{Martin Lind}
\author[1]{Adrian Muntean\thanks{adrian.muntean@kau.se}}
\author[1]{Omar Richardson}
\author[1]{Yosief Wondmagegne}
\affil[1]{Department of Mathematics and Computer Science, Karlstad University, Sweden}
\affil[2]{Capital City Kiel, Transportation Infrastructures Department, Germany}
\begin{document}

\maketitle

\begin{abstract}
    We develop and implement a numerical model to simulate the effect of photovoltaic asphalt on reducing the concentration of nitrogen monoxide (\ce{NO}) due to the presence of heavy traffic in an urban environment.
    The contributions in this paper are threefold: we model and simulate the spread and breakdown of pollution in an urban environment, we provide a parameter estimation process that can be used to find missing parameters, and finally, we train and compare this simulation with different data sets.
    We analyze the results and provide an outlook on further research.
\end{abstract}


\section{Introduction}
Pollution in urban environments has been a major issue for several decades, and efforts of combating it have spanned many areas of research.
Emerging technologies, such as those based on photo-catalysis, target the removal of vehicular nitrogen oxides (\ce{NO_x}) to mitigate the roadside air pollution problem. These techniques have proven effective and research in this area is still ongoing.
Recently, a lot has been done from the experimental perspective as well, see i.e.~\cite{WenguangFan2017A}
for experimental studies on visible–light activated photo-catalytic asphalt, or~\cite{Husken2009, Hunger2008, BALLARI201071, JKSikkema2015A} for photo-catalytic concrete products as well as~\cite{Guerrini2012165} for actual observational studies.
A comprehensive overview of the underlying principles in photo-catalysis processes in connection with removal of air pollutants is presented in~\cite{ANGELO2013522, LASEK201329, OCHIAI2012247} and the references therein.
However, from the mathematical modelling and simulation point of view, this setting is less studied.
A remotely related situation is handled via a fluid dynamics-based study on pollutant propagation in the near ground atmospheric layer as discussed in~\cite{Sukhinov2015}.
We strongly believe that there is need for a deeper insight via mathematical modelling and simulations in the connections between experiment, theory and practical applications of the methodology.
This is the place we wish to contribute. Similar techniques as the ones we apply are detailed in \cite{muntean15}, to which we refer the interested reader for more details on the modelling structure.

In this paper, we report on the use of numerical simulations to mimic the effect of
the presence of a street paved with photovoltaic asphalt has on the \ce{NO} reduction in the local ambient.
Comparisons are based on data collected at a neighborhood of one of the streets in Kiel, Germany.
The findings in this work contribute to enhance the understanding on the interplay between the various factors involved in air pollution control.
The contributions of this paper are threefold:
\begin{itemize}
    \item We model and simulate the spread and breakdown of pollution in an urban environment.
    \item We provide a parameter estimation process that can be used to determine relevant missing parameters.
    \item We train and compare this simulation with different data sets.
\end{itemize}

The paper is structured as follows: first, we present the model equations in a dimensionless formulation and explain what they are describing. Then, we provide insight in the reference set of parameters we use. As a next step, we use the designed model to formulate and solve a parameter identification problem required to identify the effect of the local environment on the \no evolution. Finally, we compare our numerical results with measured data from both before and after the photovoltaic asphalt was placed.


\section{NO pollution model}\label{model}
\label{sec:model}
\subsection{Model description}
The urban environment we model is the cross-section of a road.
This is a three-lane auto-way within city limits. In our two-dimensional model, we represent the introduction of pollution by the cars, which is diffused through the air. We model sunlight as a reactive factor~\cite{ANGELO2013522, LASEK201329}.
Finally, we model interaction with the rest of the environment by using Robin boundary conditions. These boundary conditions include environmental effects from neighbouring locations (especially molecular diffusion and dispersion), represented in a single parameter later referred to as $\sigma$. We do not incorporate wind flux in our model. Our motivation for this is two-fold: (1) we are not aware of measurements of wind data in the region under observation; (2) our model can handle a slow-to-moderate wind by a simple translation of the fluxes with an averaged convective term. However, the size effects that a strong drift introduces would be too significant for our model to capture in the current setting and with the current data available. Introducing genuine wind effects in our model would necessarily force us to handle a big area of the urban environment, where more localized measurements would be needed to determine a correct \no pollution level. Instead we opt for introducing some level of dispersion correcting the  molecular diffusion of the main pollutant agent; see Section \ref{disc} for a further discussion of this topic.

\subsection{Setting of the model equations}
\begin{figure}[ht]
\centering
\includegraphics[width=0.7\textwidth]{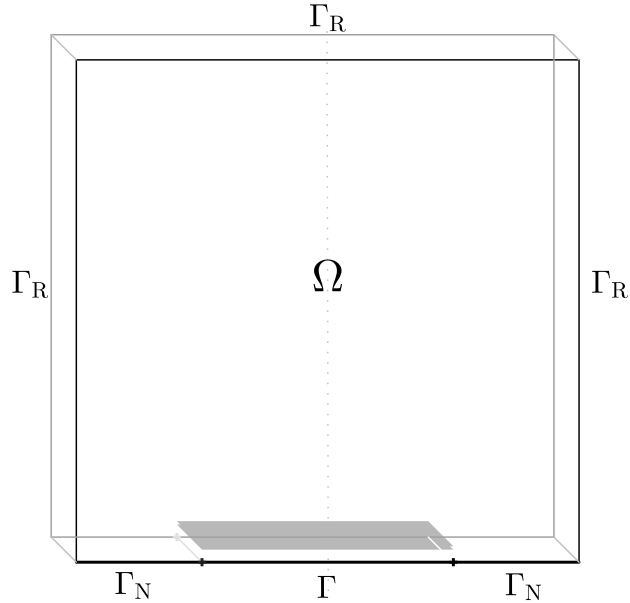}
\caption{Sketch of the cross-section. $\Gamma$ represents the face of the asphalt, the box is a cross-section of the street.
The grey rectangular box inside of the cross-section illustrates the location of \no emission.}
\label{figure:FigB1}
\end{figure}

The equations are defined in the aforementioned urban environment, denoted by $\Omega$.
The geometric representation of the environment is illustrated in Figure~\ref{figure:FigB1}.
To keep the presentation concise, we introduce the model directly into dimensionless form. We refer the reader to the Appendix for a concise explanation of the non-dimensionalization procedure.

Let $\vec{x}$ be the variable denoting the position in space and $t \in [0,1)$ be the variable denoting the time of day.
The unknown concentration profile is denoted by $u(\vec{x}, t)$ and represents the \no concentration at position $\vec{x}$ and time $t$. We refer here to \no as \emph{air pollutant}. We choose $L$ to represent a reference length scale: the width of the two times three-lane highway plus the two corresponding banks (see Section~\ref{sec:parameters}). $t_r$ is a reference time scale.

Concerning the concentration of the air pollutant, let
$u_0$ denote the initial concentration value,
$u_T$ a preset (threshold) value and $u_r$ a reference concentration value.
The preset concentration of the pollutant might correspond, for example, to the concentration of \no naturally present in the environment.

The transport (dispersion) coefficient for \no is denoted by $D$, while $T=0$ and $T=1$ denote the dimensionless initial and final times of the process observation; the start and the end of one day, respectively.
Let $f$ denote the traffic intensity on the street, measuring the average number of vehicles passing through the observation point during a specified period, with $f_r$ as a reference value.
Let $s$ denote the effect of solar radiation, with $s_r$ as a reference value. \
Introducing non-dimensional variables and rewriting gives the main equation of this study:
\begin{equation}
    \frac{\partial  u}{\partial  t} - \nabla \cdot \big(  \nabla  u  \big)  = A_f  f (\vec{x}, t) - \kappa  A_s s(t)  u,    \label{MEq1}
\end{equation}
where $A_f \, f$ represents the contribution from the emission of \no by motor vehicles.
The dependence of the coefficient $A_f$ on the other parameters can be expressed as $A_f = \frac{f_r\,L^2}{u_r \, {D}}=5.5$, a Damk\"ohler-like number that expresses the relation between the \no emitted and the density of the traffic.
The baseline value for the molecular diffusion coefficient $D$ is $0.146 \textrm{ cm}^2/\textrm{s}$; this reference value is taken from \cite{AEA}. To account for the effects of  dispersion and slow winds, we use values for this  coefficient that are of order $10^2$ higher. This brings the numerical output in the range of emission measurements for large ranges of all the other model parameters.
To quantify the evolution of the air pollutant, one has first to get a grip on the typical sizes of $\sigma$, $\kappa$ and $\gamma$.
The term $A_s\, s(t)$ represents the contribution from the reaction where \no is converted to \noo with reaction rate $\kappa$.
The dependence of $A_s$ on the other parameters is given by
$A_s = \frac{s_r\,L^2}{D}$.

The following initial and boundary conditions are imposed on \eqref{MEq1}.
\begin{subequations}
\begin{align}
    u(\vec{x}  , 0) &= \frac{u_0}{u_r} \ \
    \textrm{ for } \ \vec{x} \in \Omega, \label{MEqCond1} \\
    -\nabla u \cdot \vec{n} &= 0  \ \
    \textrm{ at } \ \Gamma_{N} \times (0  ,  1), \label{MEqCond2} \\
    \nabla u \cdot \vec{n} &= \frac{\sigma \, L}{{D}}  \left( u - \frac{u_T}{u_r} \right)^{+}
    \textrm{ at } \ \Gamma_{R} \times (0  ,  1), \label{MEqCond3} \\
    \nabla u \cdot \vec{n} &= \frac{\gamma \, \kappa \, L}{D}  u \ \
    \textrm{ at } \ \Gamma \times (0  ,  1),
    \label{MEqCond4}
\end{align}
\end{subequations}
where $g^+$ denotes the positive part of $g$, namely
$g^+(\vec{x}) = \max(g(\vec{x}),0)$ for $\vec{x} \in \Gamma_{R}$.

We refer to $\sigma$ as the \emph{environmental parameter}, to $\kappa$ as the \emph{reaction rate} and to $\gamma$ as the \emph{asphalt reactivity}. All parameters can be seen as mass-transfer coefficients; $\sigma$ represents the exchange of \no with the ambient atmosphere, $\kappa$ expresses the speed of the reaction from \no to \noo while $\gamma$ expresses the capacity of the photo-voltaic asphalt. In our model, we choose a value of $\sigma=300\textrm{ m}^3/ \mu\textrm{g}$, consistent with the base level of \no concentration in the environment according to the measurements. When simulating the scenario prior to the photo-voltaic asphalt, $\gamma=0$.

It is worth mentioning that $\kappa$ and $\gamma$ are influenced by a multitude of effects, including but not limited to the local atmospheric conditions, effect of UV, temperature and humidity as well as on the porosity and the chemical composition of the asphalt. For this reason, these parameters are situation specific and require tuning for each scenario that one wants to model.
We will do so by applying a parameter identification technique, described in Section~\ref{sec:param_ident}.
\subsection{Other parameters}
\label{sec:parameters}
\begin{itemize}
    \item Initial and threshold concentrations of \no

    $u_0$ represents the initial mass concentration of \no present in the environment. In this simulation, we choose a concentration of $u_0 = 37$ $\mu g/m^\textrm{3}$. This value corresponds with the lowest available \no concentration level from the measurements.
    $u_T$ represents the threshold concentration level from which \no disperses out of the environment we consider. In this model, we choose $u_T=0\ \mu g/m^\textrm{3}$, which means we assume low ambient levels of \no, ensuring natural dispersion even at low concentrations. In the case of high ambient levels of \no due to e.g.~nearby factories or other highways, $u_T$ can be higher.

    \item Patch-wise vehicle distribution

    The emission of \no is proportional to the amount of motor vehicles passing through this cross-section. The distribution of motor vehicles is derived from measurements of \cite{bast.de}, a German municipal service that performs automated traffic counts for a large number of cities. The traffic count we collect our data from is located at 3.5 kilometers from the \no measuring point. Because the roads in question are similar in size and geographically close, we expect this data to provide a reliable estimate. Daily, 72,278 cars passed were counted in both directions of the road. The traffic count is aggregated on an hourly level, so in order to obtain a vehicle distribution, we interpolate the hourly data with cubic splines and normalize the resulting function. We end up with a nominal density $m(t)$ of cars for each time $t \in [0,1)$ such that $\int_{0}^1 m(t) dt=1$. A plot of $m(t)$ is presented in Figure~\ref{figure:Mt}.

    \begin{figure}[ht]
    \centering
    \includegraphics[width=0.7\textwidth]{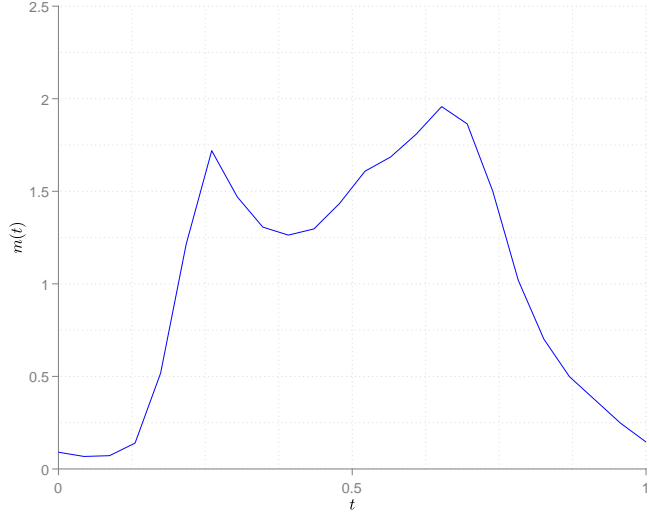}
    \caption{Traffic density $m(t)$ as a function of $t$.}
    \label{figure:Mt}
    \end{figure}

    \item Patch-wise \no emission (per vehicle)

    We model the emission of \no by choosing a specific shape for the source term $f(\vec{x},t)$ from \eqref{MEq1}. This source term is defined as follows:
    \begin{equation}
        f(\vec{x},t) = \begin{cases}
        m(t), \ \mbox{ if } \vec{x} \in A \\
        0 ,   \ \mbox{ if } \vec{x} \notin A \end{cases},
    \end{equation}
    for all $\vec{x}\in\Omega$ and $t\in[0,1)$, where $A$ is a rectangle of $15 \times 0.4$ square meters, located 0.1 meter over the asphalt. This box represents the location of emission of the vehicles. It is illustrated by the grey box in Figure~\ref{figure:FigB1}.

   \item Cross-section geometry

    The dimensions of the simulated cross-section (displayed in Figure~\ref{figure:FigB1}) are 40 meters by 8 meters.
    In the simulation, the road has a width of 15 meters ($2 \times 3$ lanes of each 2.5 meters wide) and is located in the middle of the cross section.
    This is a simplified representation of the road under consideration, but can be generalized to model roads of any dimension.
    The measuring point we use to evaluate the simulated  \no concentration profile is located 1.75 meters above the center of the road, in accordance to the real measuring point that provided us the \no concentration data.

    \item Solar effects on \no

    The natural conversion from \no to \noo due to sunlight is represented by term $\kappa A_s s(t)$. In this term, the intensity of the UV radiation is expressed by the dimensionless factor $s(t)$. We compute this factor by interpolating sunrise, sunset and solar noon data from \cite{nasa} and normalizing to obtain a function $s(t)$ such that $\int_0^1s(t)dt=1$, assuming no UV radiation outside of sunrise and sunset, and maximum UV radiation at solar noon.
    For our simulation, we choose a reference value $s_r = 1\ \textrm{UVI}$ (see Appendix), where by $\textrm{UVI}$ we mean the UV Index of the sunlight, a dimensionless quantity.
    A plot of the shape of $s(t)$ is presented in Figure~\ref{figure:St}.

\begin{figure}[ht]
\centering
\includegraphics[width=0.7\textwidth]{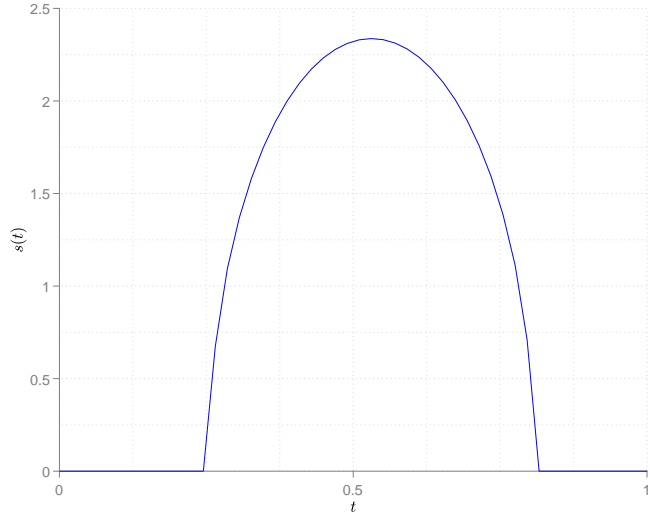}
\caption{Nondimensional UV strength $s(t)$ as a function of $t$.}
\label{figure:St}
\end{figure}

    \item Pollution-reducing effects

    As mentioned, the values of $\kappa$ and $\gamma$ strongly depend on the environment. By using the measurement data detailed in the next section, we have enough information to derive the values of these parameters for our model.
\end{itemize}

\section{Measurements}
\label{sec:measurements}
The \no measurements used in the simulation cover a period from 2012 to 2017, where \no was measured in 30-minute intervals.
We clean the data by limiting ourselves to a specific period: September 1st until December 10th: 101 measurements outside of the holiday period, with an average intensity of the sun (between summer and winter).
We interpolate this data to obtain the emission profile of a averaged day.
Data in this period from 2017 corresponds to the newly installed photovoltaic asphalt.
To compare this scenario with the initial situation, we use measurement data from 2016.

\begin{figure}[ht]
\centering
\includegraphics[width=0.7\textwidth]{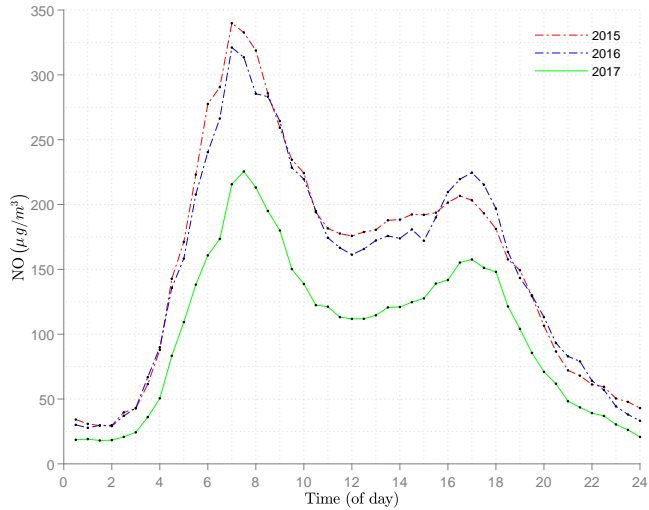}
\caption{\no measurements at the Theodorus-Heuss-Ring in Kiel, Germany. The green curve (2017) has the same shape as the other two curves, but presents a significant reduction of \no levels.}
\label{figure:FigA2}
\end{figure}

To assert that the traffic intensity remained relatively constant, we compare the measurements from 2016 to the same period in 2015.
Figure~\ref{figure:FigA2} shows comparable concentration profiles.
This indicates, as a reasonable assumption, that the only major change between 2016 and 2017 was the new asphalt, given the local weather conditions incorporated in the parameter $\kappa$.
Figure~\ref{figure:FigA2} reveals a large reduction in \no pollution after 2016.

\section{Parameter identification}
\label{sec:param_ident}
As stated in Section~\ref{sec:model}, we require situation-specific values for the parameters $\kappa$ and $\gamma$. To obtain these, we propose and solve a \emph{parameter identification problem} based on the datasets described in Section~\ref{sec:measurements}.
In a two-step procedure, we first use the measurements prior to the installation to obtain $\kappa$, knowing that in this case $\gamma=0$, and then use the measurements after the installation to obtain $\gamma$.

Mathematically speaking, we wish to solve the following problem:
Let $u(\cdot ; \kappa , \gamma)$ be the solution to \eqref{MEq1} for given parameters $\kappa$ and $\gamma$ and $u_r$ be the measurement data.
Find $\kappa$ and $\gamma$ that solve the following optimization problem:
\begin{equation}
    \min_{\kappa \in \mathcal{P}_{\kappa}}
    \min_{\gamma \in \mathcal{P}_{\gamma}}
    \|u_r(\vec{x},t) - u(\vec{x},t;\kappa,\gamma) \|_{\mathcal{L}^2(\Omega \times (0,1))}.
\end{equation}
Here $\mathcal{P}_{\kappa}$ and $\mathcal{P}_{\gamma}$ are compact sets in $\mathbb{R}$ where parameters are searched.

\section{Simulation}
\label{sec:simulation}
This section describes the setup of determining the effectiveness of the photovoltaic asphalt.

\subsection{Simulation framework}

We numerically compute the solution to \eqref{MEq1} with \eqref{MEqCond1}-\eqref{MEqCond4} with a finite element simulation using the FEniCS library~\cite{fenics}. We investigate the process within the specified cross-section, leading to a two-dimensional reaction-diffusion process. The finite element mesh has $30\times 30$ elements on a rectangular grid. The solution to \eqref{MEq1} is approximated with quadratic basis functions.

We simulate two scenarios:
(A) the \no concentration profile prior to the photovoltaic asphalt (corresponding to measurements from 2015 and 2016); and
(B) the \no concentration profile after the installation of the photovoltaic asphalt (corresponding to measurements from 2017).
In these simulations, we use the parameters as described in Section~\ref{sec:model} and choose a diffusion coefficient of $D=43.8$, which shows agreement with the current setting. In case (A) specifically, we fix $\gamma=0$.

Our goal is to use the data set from (A) to train our simulation on the value of parameter $\kappa$. Then we use the data set from (B) to determine the effect the photovoltaic asphalt has on the reduction of the local \no concentration, i.e., find $\gamma$.

\begin{figure}[!ht]
\centering
\includegraphics[width=0.7\textwidth]{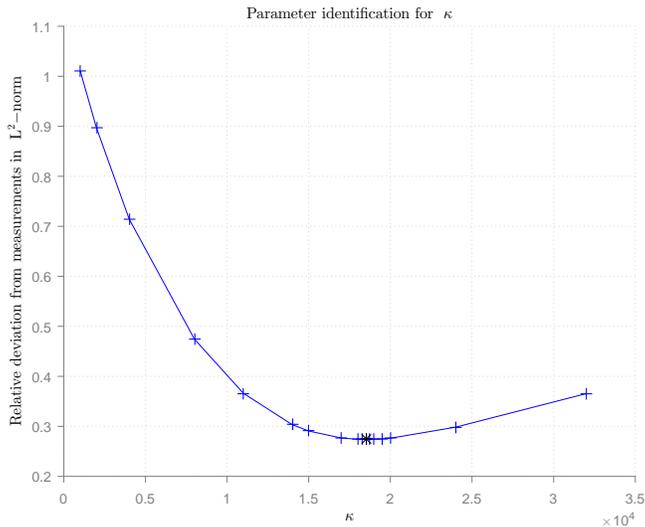}
\caption{The relative discrepancy between the simulated \no concentrations and the measurements for a range of $\kappa$ values and $\gamma=0$. The optimum is reached for $\kappa = 1.85\times 10^4\ 1/[\textrm{day}\ \textrm{UVI}]$, displayed as the minimum of this curve.}
\label{figure:opt_kappa}
\end{figure}

Figure~\ref{figure:opt_kappa} displays this process, where for a range of $\kappa$ we plotted the relative discrepancy with the measurements.
More precisely, let $u_r$ denote the interpolation of the measurements and $u_\kappa$ the result of our simulation for a specific $\kappa$, then the relative discrepancy $e_\kappa(u)$ is defined as

\begin{equation}
    e_\kappa(u) = \frac{||u_r-u_\kappa||_{\mathcal{L}^2(\Omega \times (0,1))}}{||u_r||_{\mathcal{L}^2(\Omega \times (0,1))}}.
    \label{eq:rel_dev}
\end{equation}

Having captured the nature of the process, we can estimate the effect of the photovoltaic asphalt in terms of $\gamma$ by starting a new parameter identification process using the data set from the second scenario.

\begin{figure}[!ht]
\centering
\includegraphics[width=0.7\textwidth]{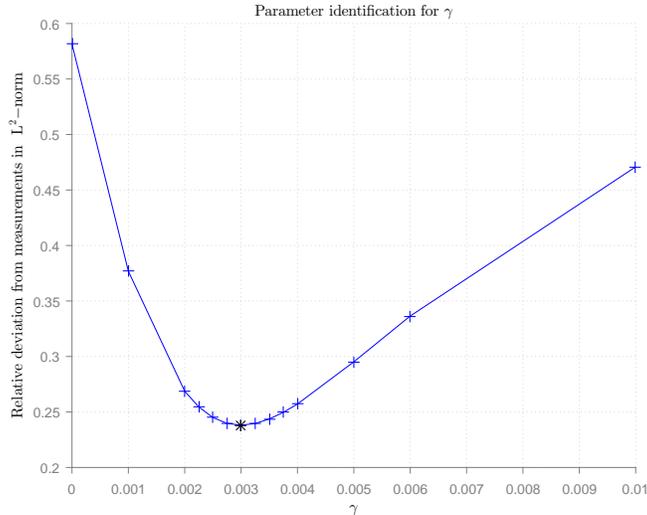}
\caption{The relative discrepancy between the simulated \no concentrations and the measurements for a range of $\gamma$ values and $\kappa = 1.85\times 10^4\ 1/[[\textrm{day}\ \textrm{UVI}]$. The optimum is reached for $\gamma = 3.0\times 10^{-3}$, displayed as the minimum of this curve.}
\label{figure:opt_gamma}
\end{figure}

Figure~\ref{figure:opt_gamma} shows the relative discrepancy (similarly defined as in \eqref{eq:rel_dev}) for $u_\gamma$ (with an optimal value of $\kappa$).
Figure~\ref{figure:opt_kappa} and Figure~\ref{figure:opt_gamma} suggest that we deal with a convex minimization problem.

\begin{figure}[ht]
\centering
\includegraphics[width=0.7\textwidth]{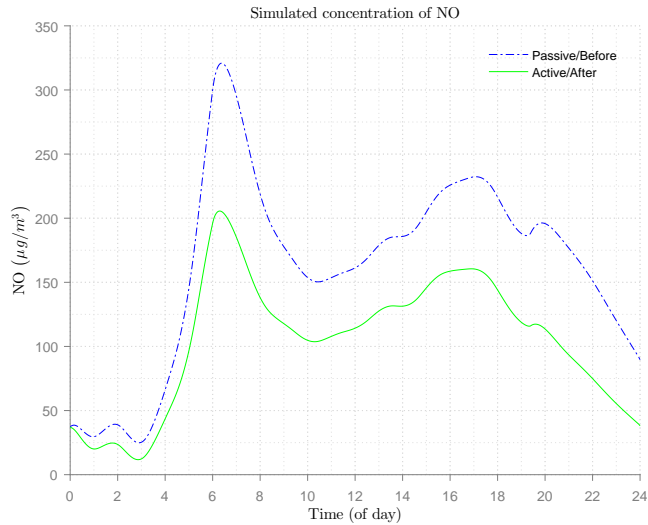}
\caption{Simulated \no concentration profiles. The dash-dotted line corresponds to the period in 2016, the continuous line has corresponds to the period in 2017.}
\label{figure:sim_no}
\end{figure}

The simulations with optimally fitted parameters, before and after the construction of the photovoltaic asphalt, are displayed in Figure~\ref{figure:sim_no}.

\begin{figure}[ht]
\centering
\includegraphics[width=0.7\textwidth]{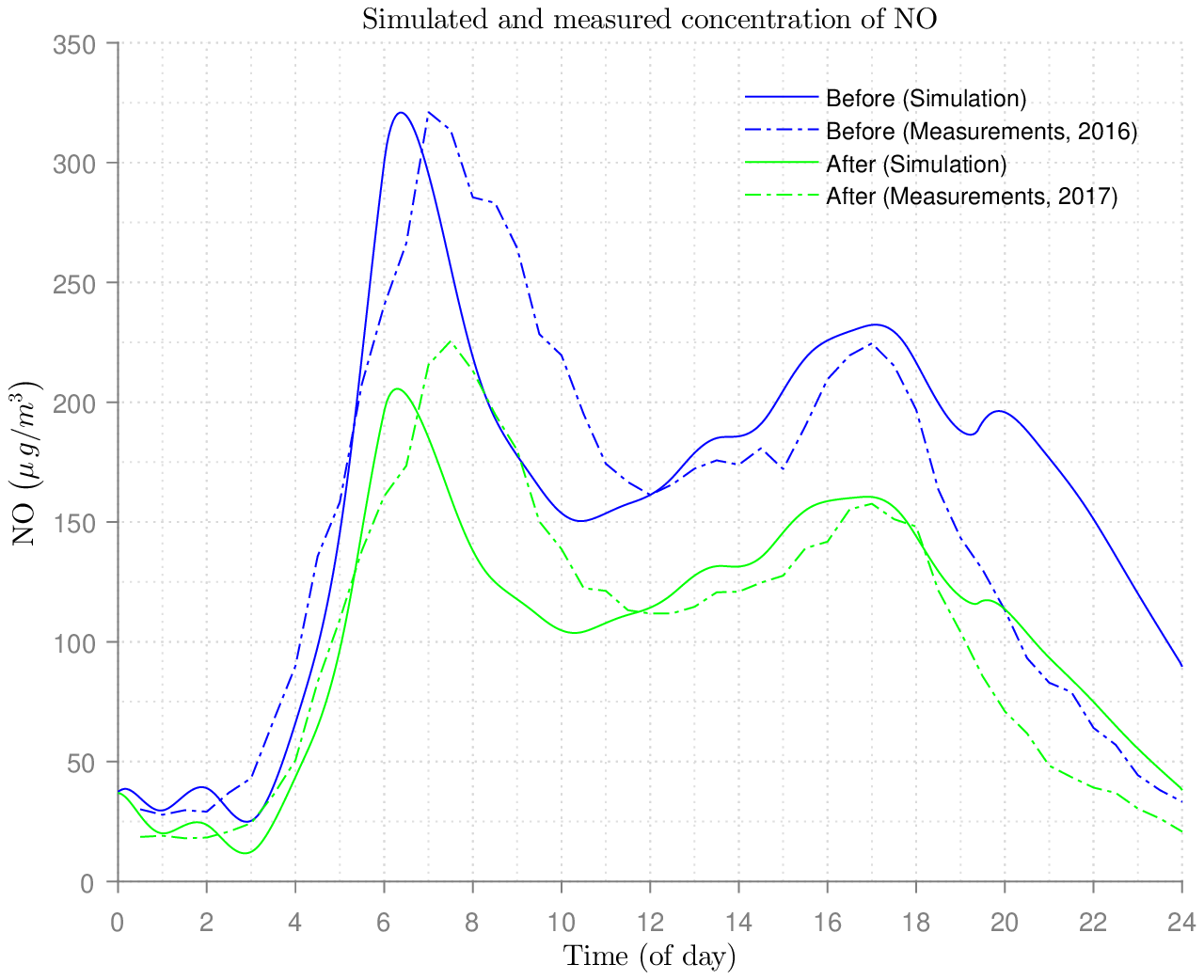}
\caption{\no measurements compared to simulations in the two cases. The simulations show good quantitative agreement.}
\label{fig:comparison}
\end{figure}

We use the measurements of \no to compare the simulation results with. Figure~\ref{fig:comparison} shows the measured and simulated data combined. The simulation agrees qualitatively and quantitatively, with an mass error of 58.64 $\mu g/m^3$ in the pre-photovoltaic case and 22.68 $\mu g/m^3$ in the post-photovoltaic one.
This mass error is defined as
\begin{equation}
    \int_0^T\int_\Omega \left|u(\vec{x},t) - u_r(\vec{x},t) \right|d\vec{x}dt
\end{equation}
where $u$ represents the simulation solution and $u_r$ the measurements.

\section{Discussion and outlook}\label{disc}
This paper shows that it is possible to use a mathematical model like the one presented in Section \ref{model} to describe the evolution of \no concentration as a function of time in an urban environment where both intense motorized vehicle traffic and photovoltaic asphalts are present.

A number of things can be done to continue this investigation and improve the quantitative predictions obtained in this study.
\begin{itemize}
    \item Perform a 3D simulation along the whole length of the photovoltaic asphalt.
    \item Account for the presence of uncertainties in the weather condition (especially the variation in UV radiation and the effect of precipation).
    \item Numerically identify scenarios leading to extreme \no pollution.
    \item By neglecting wind flux, we emphasize that our interest lies in quantifying the reactive part of this reaction-diffusion process. This approach allows us to identify two extreme scenarios: (1) excellent performance of the asphalt, and (2) low-response on \no pollution induced by the presence of the asphalt. Including more detailed wind effects is certainly of interest (both from a theoretical and a practical perspective), since it is known that the location under investigation (Kiel, Germany) is prone to varied weather conditions.
\end{itemize}

\section{Acknowledgements}
We gratefully thank the State Agency for Agriculture, Environment and Rural Areas Schleswig-Holstein (LLUR) for providing measurement data at the Theodor-Heuss-Ring in Kiel, Germany. The authors acknowledge the very valuable input from the side of the referees regarding the shaping of the final form of the manuscript.

\section{Appendix}
For $\vec{x} \in \Omega$ and $t \in [0 , T)$, we consider the following equation
\begin{equation}
\frac{\partial  u}{\partial  t} - \nabla \cdot \big(D \nabla u \big)  = f(\vec{x}, t) - \kappa\, s(t)\, u(\vec{x}, t).
\label{DL1}
\end{equation}
Let us also introduce the following rescalings:
$\bar{\vec{x}}=\frac{\vec{x}}{L}$,
$\bar{t}=\frac{t}{t_r}$,
$\bar{u}=\frac{u}{u_r}$,
$\bar{f}=\frac{f}{f_r}$ and
$\bar{s}=\frac{s}{s_r}$, where
\[f_r = \max_{{(\vec{x} , t)\in\Omega\times [0,T]}}|f(\vec{x}, t)| \ . \]
Rewriting \eqref{DL1} using these new rescaled functions and variables yields

\begin{equation}
\frac{\partial \bar{u}}{\partial \bar{t}}  -
\frac{t_r D}{L^2} \bar{\nabla} \cdot \big(\bar{\nabla} \bar{u} \big) =
\frac{t_r \, f_r}{u_r} \bar{f} - \kappa \, t_r \, s_r \, \bar{s} \, \bar{u}.
\label{DL2}
\end{equation}
This equation is posed in a rescaled space domain $\bar{\Omega}=\frac{\Omega}{L}$.
Choosing $t_r = \frac{L^2}{D}$ reduces \eqref{DL2} to the form
\begin{equation}
\frac{\partial \bar{u}}{\partial \bar{t}}  -
\bar{\nabla} \cdot \big(\bar{\nabla} \bar{u} \big) =
\frac{L^2 \, f_r}{D \, u_r} \bar{f} - \kappa \, \frac{L^2 \, s_r }{D} \bar{s} \, \bar{u}.
\label{DL3}
\end{equation}
Setting $A_f=\frac{L^2 \, f_r}{D \, u_r}$, $A_s = \frac{L^2 \, s_r }{D}$ and abusing notation, namely rewriting $\Omega$ instead of $\bar{\Omega}$, $u$ instead of $\bar{u}$,
$f$ instead of $\bar{f}$, and $s$ for $\bar{s}$ in \eqref{DL3} posed in $\bar{\Omega}$
gives \eqref{MEq1} posed in $\Omega$. A similar discussion yields the structure of the initial and boundary conditions in non-dimensional form.

\bibliographystyle{actapoly}
\bibliography{NOxBibl}

\begin{thebibliography}{10}

\bibitem{fenics}
Martin Alnæs, Jan Blechta, Johan Hake, August Johansson, Benjamin Kehlet,
  Anders Logg, Chris Richardson, Johannes Ring, Marie Rognes, and Garth Wells.
\newblock The {FEniCS} {P}roject {V}ersion 1.5.
\newblock {\em Archive of Numerical Software}, 3(100), 2015.

\bibitem{BALLARI201071}
M.M. Ballari, M.~Hunger, G.~Hüsken, and H.~Brouwers.
\newblock Modelling and experimental study of the $\textrm{NOx}$ photocatalytic
  degradation employing concrete pavement with titanium dioxide.
\newblock {\em Catalysis Today}, 151(1):71 -- 76, 2010.
\newblock 2nd European Conference on Environmental Applications of Advanced
  Oxidation Processes(EAAOP-2).

\bibitem{WenguangFan2017A}
Wenguang Fan, Ka~Yan Chan, Chengxu Zhang, and Michael~K.H. Leung.
\newblock {Advanced solar photocatalytic asphalt for removal of vehicular NOx}.
\newblock {\em Energy Procedia}, 143:811 -- 816, 2017.

\bibitem{bast.de}
{Federal Ministry of Transport and Digital Infrastructure}.
\newblock Federal {H}ighway {R}esearch {I}nstitute {(BASt)}.
\newblock \url{http://www.bast.de/DE/Home/home_node.html}.
\newblock Accessed: 2018-01-16.

\bibitem{nasa}
{Global Radiation Group}.
\newblock {Earth System Research Laboratory; Global Monitoring Division}.
\newblock \url {https://www.esrl.noaa.gov/gmd/grad/solcalc/}.
\newblock Accessed: 2018-01-16.

\bibitem{Guerrini2012165}
Gian~Luca Guerrini.
\newblock {Photocatalytic performances in a city tunnel in Rome: NOx monitoring
  results}.
\newblock {\em Construction and Building Materials}, 27(1):165 -- 175, 2012.

\bibitem{Hunger2008}
M.~Hunger, H.J.H. Brouwers, and M.~Ballari.
\newblock Photocatalytic degradation ability of cementitious materials: A
  modeling approach.
\newblock In W.~Sun, K.~{Breugel, van}, C.~Miao, G.~Ye, and H.~Chen, editors,
  {\em Proceedings of the 1st International Conference on Microstructure
  related Durability of Cementitious Composites, 13 -15 October 2008, Nanjing,
  China.}, pages 1103--1112, 2008.

\bibitem{Husken2009}
G.~H{\"u}sken, M.~Hunger, and H.~Brouwers.
\newblock Experimental study of photocatalytic concrete products for air
  purification.
\newblock {\em Building and environment}, 44(12):2463--2474, 2009.

\bibitem{LASEK201329}
Janusz Lasek, Yi-Hui Yu, and Jeffrey~C.S. Wu.
\newblock {Removal of $\textrm{NOx}$ by photocatalytic processes}.
\newblock {\em Journal of Photochemistry and Photobiology C: Photochemistry
  Reviews}, 14:29 -- 52, 2013.

\bibitem{muntean15}
Adrian Muntean.
\newblock {\em Continuum Modeling: An Approach Through Practical Examples}.
\newblock Springer, 2015.

\bibitem{OCHIAI2012247}
Tsuyoshi Ochiai and Akira Fujishima.
\newblock {Photoelectrochemical properties of $\textrm{TiO}_2$ photocatalyst
  and its applications for environmental purification}.
\newblock {\em Journal of Photochemistry and Photobiology C: Photochemistry
  Reviews}, 13(4):247 -- 262, 2012.

\bibitem{JKSikkema2015A}
J.K. Sikkema, S.K. Ong, and J.E. Alleman.
\newblock {Photocatalytic concrete pavements: laboratory investigation of
  $\textrm{NO}$ oxidation rate under varied environmental conditions}.
\newblock {\em Construction and Building Materials}, 100:305 -- 314, 2015.

\bibitem{Sukhinov2015}
A.~I. Sukhinov, D.~S. Khachunts, and A.~E. Chistyakov.
\newblock A mathematical model of pollutant propagation in near-ground
  atmospheric layer of a coastal region and its software implementation.
\newblock {\em Computational Mathematics and Mathematical Physics},
  55(7):1216--1231, Jul 2015.

\bibitem{AEA}
J.~Targa and A.~Loader.
\newblock Diffusion tubes for ambient {$\textrm{NO2}$} monitoring: Practical
  guidance.
\newblock Technical report, AEA Energy \& Environment, 2008.

\bibitem{ANGELO2013522}
Joana Ângelo, Luísa Andrade, Luís~M. Madeira, and Adélio Mendes.
\newblock {An overview of photocatalysis phenomena applied to $\textrm{NOx}$
  abatement}.
\newblock {\em Journal of Environmental Management}, 129:522 -- 539, 2013.

\end{thebibliography}

\end{document}